\documentclass[headsepline,dvips]{scrartcl}
\usepackage{latexsym, amssymb,amsmath}
\def\SL2{\operatorname{SL}_{2}(K)}

\def\GL2{\operatorname{GL}_{2}(K)}
\def\Ga{{\mathbb G}_{a}}

\def\INVSL2{$K[V]^{operatorname{SL}_{2}(K)}$}
\def\INVSO2{$K[V]^{operatorname{SO}_{2}(K)}$}
\def\INVGL2{$K[V]^{operatorname{GL}_{2}(K)}$}

\def\Bew{\noindent \textit{Proof. }}
\def\qed{\hfill $\Box$}
\def\Magma{{\sc Magma }} 

\def\depth{\operatorname{depth}}
\def\height{\operatorname{height}}
\def\cmdef{\operatorname{cmdef}}

\def\Hom{\operatorname{Hom}}
\def\GL{\operatorname{GL}}
\def\SL{\operatorname{SL}}
\def\id{\operatorname{id}}

\def\Ass{\operatorname{Ass}}
\def\Ann{\operatorname{Ann}}

\def\myforall{\textrm{ for all }}

\newtheorem{Lemma}{Lemma}[section]
\newtheorem{Satz}[Lemma]{Theorem}
\newtheorem{Def}[Lemma]{Definition}
\newtheorem{Korollar}[Lemma]{Corollary}

\newtheorem{Prop}[Lemma]{Proposition}

\newtheorem{Bsp}[Lemma]{Example}

\title{On the depth of invariant rings of infinite groups}
\author{Martin Kohls \\
Technische Universit\"at M\"unchen, 
 Zentrum Mathematik - M11\\
Boltzmannstrasse 3, 
 85748 Garching, Germany\\
\\
kohls@ma.tum.de}

\begin{document}
\maketitle
\begin{abstract}
Let $K$ be an algebraically closed field. For a finitely generated graded
commutative $K$-algebra $R$, let
$\cmdef R:=\dim R-\depth R$ denote the Cohen-Macaulay defect of $R$. Let $G$ be a linear algebraic group
over $K$ that is reductive but not linearly reductive. We show that there
exists a faithful rational representation $V$ of $G$ (which we will give explicitly) such that $\cmdef
K\left[V^{\oplus k}\right]^{G}\ge k-2$ for all $k\in {\mathbb N}$. 
\end{abstract}

\section{Introduction}

Let us first fix some notation and conventions. The symbol $G$ will always denote a linear algebraic group over the
algebraically closed field $K$ of characteristic $p\ge 0$, and $V$ will always
denote a finite dimensional rational $G$-module.  Then $G$ acts on the ring $K[V]$ of polynomial functions $V\rightarrow K$ by $\sigma\cdot f := f\circ \sigma^{-1}$,
where $f\in K[V],\sigma\in G$, so we have $K[V]\cong S(V^{*})$, the symmetric
algebra of the dual $V^{*}=\Hom_{K}(V,K)$. We will always assume that
the algebra $K[V]^{G}$ of invariants is finitely generated, which due to Nagata
\cite{NagataAffine} is always true if $G$ is reductive. We are interested in how far
the invariant ring $K[V]^{G}$ is from being Cohen-Macaulay. To this end we
define the Cohen-Macaulay defect $\cmdef K[V]^{G}$ to be $\dim
K[V]^{G}-\depth K[V]^{G}$. Then $\cmdef K[V]^{G}\ge 0$ and $K[V]^{G}$ is
Cohen-Macaulay precisely when $\cmdef K[V]^{G}=0$. By Hochster and Roberts
\cite{HochsterRoberts}, the invariant ring $K[V]^{G}$ is Cohen-Macaulay if $G$
is linearly reductive. 

Several papers have dealt with the depth of invariant
rings of \emph{finite} groups, see for example  
\cite{CampEtAl,Ellingsrud,FleischmannCohomology,FleischmannCohomConnectivity,FleischmannDepthTrace,Gordeev,KemperDepthCohom,NeuselDepth,ShankWehlau}, but there are up to today no quantitative results for the
depth of invariants of (infinite) algebraic groups.

 The goal of this paper is to make a
first step in this direction. To have a context for the results of this paper,
we mention three results from the literature explicitly:
Every overview has to start with the celebrated theorem of Ellingsrud and Skjelbred \cite{Ellingsrud}, who proved
that for $p>0$ and $G$ a cyclic $p$-group, we have
$\cmdef K[V]^{G}=\max\left(\dim_{K} (V)-\dim_{K}(V^{G})-2, 0\right)$. The main part of
our paper is about bringing the following two results together:

\begin{Satz}[{Gordeev, Kemper \cite[Corollary 5.15]{Gordeev}}] \label{GordiKemper}
Let $G$ be finite and $p||G|$. Then for every faithful $G$-module $V$, we have
\[
\lim_{k\to \infty}\cmdef K\left[V^{\oplus k}\right]^{G}=\infty.
\]
\end{Satz}
(Where $V^{\oplus k}$ denotes the $k$-fold direct sum of $V$.)

\begin{Satz}[{Kemper \cite[Theorem 7]{KemperLinRed}}]
Let  $G$ be reductive, but not linearly reductive. Then there exists a $G$-module
$V$ such that $K[V]^{G}$ is not Cohen-Macaulay. 
\end{Satz}

These two theorems should be compared with the following theorem, which is the
main result of our paper:
\begin{Satz}\label{DasHauptErgebnis}
Let  $G$ be reductive, but not linearly reductive.
 Then there exists a faithful $G$-module
$V$ (that will be given explicitly in the proof) such that 
\[
\cmdef K\left[V^{\oplus k}\right]^{G}\ge k-2 \quad \myforall k\ge 1.
\]
\end{Satz}

On the other hand, if $G=\SL_{n}$ acts on $V=K^{n}$ by left multiplication,
then $V$ is a faithful $\SL_{n}$-module, and if $n\ge 2$ and $p>0$, then
$\SL_{n}$ is not linearly reductive. But by Hochster \cite[Corollary 3.2]{Hochster}, we
have $\cmdef K\left[V^{\oplus k}\right]^{\SL_{n}}=0$ for all $k\in{\mathbb
  N}$, so Theorem \ref{GordiKemper} can not be generalized to reductive groups that are
not linearly reductive.

It would be interesting to know if $\cmdef K[V]^{G}> 0$ implies $\lim_{k\to \infty}\cmdef K\left[V^{\oplus
    k}\right]^{G}=\infty$, but we have
neither a proof nor a counterexample for this.\\

{\bf Acknowledgment.} This paper contains the main result of my
``Doktorarbeit'' (Ph.D. thesis) \cite{Doktorarbeit} (in German). I want to take the opportunity to thank my
advisor Gregor Kemper for his constant support.

\section{The depth of graded algebras}
For the convenience of the reader, we have collected some standard facts about
the depth of graded algebras that can be looked up in any better book on
commutative algebra like \cite{BrunsHerzog,Eisenbud}. The given references
generally only treat the local case, but this case carries over to our graded
situation. For this paper, $R$ will denote a finitely generated graded
commutative $K$-algebra $R=\bigoplus_{d=0}^{\infty}R_{d}$ where $R_{0}=K$.
We call $R_{+}:=\bigoplus_{d=1}^{\infty}R_{d}$ the \emph{maximal
  homogeneous ideal of R}. 
A sequence of homogeneous elements
$a_{1},\ldots,a_{k}\in R_{+}$ is called a \emph{partial homogeneous system of
  parameters (phsop)} if $\height (a_{1},\ldots,a_{k})_{R}=k$. If $k=\dim R$,
then the sequence is called a \emph{homogeneous system of
  parameters (hsop)}, and then $R$ is finitely generated as a module over
$A:=K[a_{1},\ldots,a_{k}]$.  Due to the Noether normalization theorem, hsops
always exist.

The following lemma makes it easier to find phsops in an invariant ring.

\begin{Lemma}[{Kemper \cite[Lemma 4]{KemperLinRed}}] \label{redphsop}
Let $G$ be reductive. If $a_{1},\ldots,a_{k} \in K[V]^{G}$ form a phsop in $K[V]$, then they also
form a phsop in $K[V]^{G}$.
\end{Lemma}

Let $M$ always denote a nonzero, finitely generated ${\mathbb Z}$-graded
$R$-module (the most important case is $M=R$).
Then a sequence of homogeneous elements
$a_{1},\ldots,a_{k}\in R_{+}$ is called a \emph{homogeneous $M$-regular
  sequence of length k} if for each
$i=1,\ldots,k$ we have that $a_{i}$ is not a zero divisor of
$M/(a_{1},\ldots,a_{i-1})_{R}M$. If $I\subseteq R_{+}$ is a proper homogeneous 
ideal, then a homogeneous $M$-regular sequence $a_{1},\ldots,a_{k}\in I$ is
called maximal (in $I$), if it can not be extended to a longer $M$-regular
sequence lying in $I$. Due to the theorem of Rees, two such maximal
$M$-regular sequences have the same length (which is finite in our setup), and
we write $\depth (I,M)$ for that common length (the
\emph{$I$-depth of $M$}). We call $\depth M:=\depth (R_{+},M)$ the \emph{depth
  of $M$}. In the case of $M=R$, every regular sequence is a phsop, and $R$ is
Cohen-Macaulay if and only if every phsop is a regular sequence. We write $\depth I:=\depth (I,R)$ for
the \emph{depth of $I$}. We always have $\depth I\le \height I$, so for the
\emph{Cohen-Macaulay defect of $I$}, $\cmdef I:=\height I-\depth I$, we always
have $\cmdef I\ge 0$.

\begin{Satz}[{\cite[Exercise 1.2.23]{BrunsHerzog}}]\label{cmdefabschatzung}
Let $I\subseteq R_{+}$ be a proper homogeneous ideal. Then
\[
\depth R\le \depth I +\dim I.
\]
In particular we have $\cmdef R\ge\cmdef I$.
\end{Satz}

To apply  this theorem in order to get a good lower bound for $\cmdef R$, one has
to find ideals $I\subseteq R_{+}$ of which one knows the depth. The following
 lemma, which is inspired by  Shank and Wehlau \cite[Theorem
2.1]{ShankWehlau} is the proper tool.

\begin{Lemma} \label{SatzShankWehlau}
Let $I\subseteq R_{+}$ be a proper homogeneous ideal of
$R$, and $a_{1},\ldots,a_{k}\in I$ be a homogeneous $M$-regular
sequence.
Then $\depth (I,M)=k$ if and only if there exists a  $m\in M$ with $m\notin
(a_{1},\ldots,a_{k})M$ and $Im\subseteq (a_{1},\ldots,a_{k})M$.
\end{Lemma}

\Bew First assume that there exists  a $m\in M$ with $m\notin
(a_{1},\ldots,a_{k})M$ and $Im\subseteq (a_{1},\ldots,a_{k})M$. Then obviously
$I$ only consists of zero divisors of $M/(a_{1},\ldots,a_{k})M$, hence
$a_{1},\ldots,a_{k}$ is a maximal $M$-regular sequence in $I$, hence $\depth
(I,M)=k$.

Conversely assume $\depth (I,M)=k$. Then $a_{1},\ldots,a_{k}\in I$ is a
maximal homogeneous $M$-regular sequence. Then $I$ only consists of zero
divisors of $N:=M/(a_{1},\ldots,a_{k})M$, so by  \cite[Theorem
3.1.b]{Eisenbud} we have
$I\subseteq \bigcup_{\wp \in \Ass_{R}(N)}\wp$, and by prime avoidance \cite[Lemma
3.3]{Eisenbud} there is a $\wp\in\Ass_{R}N$ with $I\subseteq \wp$. Since $\wp$
is an associated prime ideal of $N$, there is a $n\in N\setminus\{0\}$ with
$\wp=\Ann_{R}n$, and for a $m\in M$
with $n=m+(a_{1},\ldots,a_{k})M$ we have $m\not\in(a_{1},\ldots,a_{k})M$ and $Im\subseteq
  (a_{1},\ldots,a_{k})M$. \qed\\

We will apply this lemma only in the case $k=2$, since it is difficult to
check if $k\ge 3$ elements form a regular sequence. To check if two elements
form a regular sequence, we have the following lemma.

\begin{Lemma} \label{ZweiReg}
Let $a_{1},a_{2}\in K[V]_{+}$ be homogeneous. Then the following are
equivalent:

(a) $a_{1},a_{2}$ form a phsop in $K[V]$.

(b) $a_{1},a_{2}$ form a regular sequence in $K[V]$.

(c) $a_{1},a_{2}$ are coprime in $K[V]$.\\
If one (hence all) of the above conditions is satisfied and we additionally
have $a_{1},a_{2}\in K[V]^{G}$, then $a_{1},a_{2}$ also form a regular
sequence in $K[V]^{G}$.
\end{Lemma}

\section{First cohomology of groups and the depth of invariants}
The results of this section are a quantitative extension of the
qualitative results of Kemper~\cite{KemperLinRed}. 
Let $H$ be an arbitrary group and $W$ be a $KH$-module (not
necessarily finite dimensional). A map $g: H\rightarrow W, \sigma\mapsto
g_{\sigma}$ is called a \emph{$(1)$-cocycle}, if we have $g_{\sigma\tau}=\sigma
g_{\tau}+g_{\sigma}$ for all $\sigma,\tau\in H$. The sum of two cocycles is
again a cocycle, so the set of all cocycles $Z^{1}(H,W)$ is an additive
group. For any $w\in W$, the map $H\rightarrow W$ given by $\sigma\mapsto (\sigma-1)w:=\sigma w
- w$ is also a cocycle, and we call a cocycle which is given by such a $w$ a
\emph{$(1)$-coboundary}. The set of all coboundaries $B^{1}(H,W)$ is obviously
a subgroup of $Z^{1}(H,W)$, and we write $H^{1}(H,W):=Z^{1}(H,W)/B^{1}(H,W)$
for the corresponding factor group. We call a cocycle $g$
non-trivial, if it is not a coboundary, and we will often confuse an element
$g\in Z^{1}(H,W)$ with its image (also denoted $g$) in $H^{1}(H,W)$. Thus $g$
is non-trivial if and only if $g\ne 0$ in $H^{1}(H,W)$.
If $H$ is a linear algebraic group and $W$ a rational (not necessarily finite
dimensional) $H$-module, then by $Z^{1}(H,W)$ we will always mean the
cocycles that are given by morphisms of $H$ to $W$ (this is automatic for $B^{1}$).
Now let $a\in K[V]^{G}$ be an invariant. Then for any $g\in Z^{1}(G,K[V])$, we
can define an element $ag\in Z^{1}(G,K[V])$ by $(ag)_{\sigma}:=ag_{\sigma}$ for all
$\sigma\in G$. Obviously, multiplication with $a$ gives a group homomorphism
$Z^{1}(G,K[V])\rightarrow Z^{1}(G,K[V])$, and induces a group homomorphism
$H^{1}(G,K[V])\rightarrow H^{1}(G,K[V])$. Now for any $g\in H^{1}(G,K[V])$, we
define its \emph{annihilator} as
\[
\Ann_{K[V]^G}(g):=\left\{a\in K[V]^{G}: a\cdot g=0\in H^{1}(G,K[V])\right\}\unlhd
K[V]^{G}.
\]
This ideal is proper if and only if $g\ne 0\in H^{1}(G,K[V])$.
We call a $0\ne g\in Z^{1}(G,K[V])$ \emph{homogeneous of degree $d\ge 0$}, if  $g_{\sigma}\in K[V]_{d}$ for all $\sigma\in G$. An element of
$H^{1}(G,K[V])\setminus\{0\}$ is called \emph{homogeneous of degree $d\ge 0$}, if it can be represented by a homogeneous
element of degree $d$ of $Z^{1}(G,K[V])\setminus B^{1}(G,K[V])$ (this is well defined). If $g\in H^{1}(G,K[V])$ is homogeneous, then its
annihilator $\Ann_{K[V]^G}(g)$ is also homogeneous.

The proof of the following proposition has some overlap with the one of
Kemper \cite[Proposition 6]{KemperLinRed}, but we get a sharper result here.

\begin{Prop} \label{depthKillerideal}
 Let $0\ne g\in H^{1}(G,K[V])$ be homogeneous, and assume there exist  two homogeneous elements $a_{1},a_{2}\in\Ann_{K[V]^G}(g)$  of positive degree that are coprime in $K[V]$. Then
\[
\depth(\Ann_{K[V]^G}(g))=2.
\]
\end{Prop}

\Bew Because of Lemma \ref{ZweiReg}, $a_{1},a_{2}$ form a $K[V]^{G}$-regular sequence in
$\Ann_{K[V]^G}(g)$. Since $a_{i}g=0\in H^{1}(G,K[V])$ (i=1,2), there are elements
$b_{1},b_{2}\in K[V]$ such that 
\[
a_{i}g_{\sigma}=(\sigma-1)b_{i}\quad \myforall \sigma\in G, i=1,2.
\]
Now set $m:=a_{1}b_{2}-a_{2}b_{1}\in K[V]^{G}$. We will show that $m$ fulfills
the hypotheses of Lemma~\ref{SatzShankWehlau} with $R=M=K[V]^{G}$, $I=\Ann_{K[V]^G}(g)$
and $k=2$,
i.e. $m\not\in(a_{1},a_{2})_{K[V]^{G}}$ and $m\Ann_{K[V]^G}(g)\subseteq
(a_{1},a_{2})_{K[V]^{G}}$. Then Lemma \ref{SatzShankWehlau}
yields $\depth (\Ann_{K[V]^G}(g))=2$.

Assume by way of contradiction $m\in (a_{1},a_{2})_{K[V]^{G}}$. Then there are $f_{1},f_{2} \in K[V]^{G}$ with 
\begin{equation*} 
m=a_{1}b_{2}-a_{2}b_{1}=f_{1}a_{1}+f_{2}a_{2}.
\end{equation*}
Then $a_{1}(b_{2}-f_{1})=a_{2}(f_{2}+b_{1})$, and
$a_{1},a_{2}$ being coprime in $K[V]$ yields that $a_{1}$ is a divisor of $f_{2}+b_{1}$, hence $f_{2}+b_{1}=a_{1} \cdot h$ with $h \in K[V]$. Now
\[
a_{1} \cdot (\sigma -1)h=(\sigma -1)(a_{1}h)=(\sigma -1)(f_{2}+b_{1})=(\sigma -1)b_{1}=a_{1}g_{\sigma} \qquad \myforall \sigma \in G,
\]
hence $g_{\sigma}=(\sigma -1)h$ for all $\sigma\in G$. But then $g=0\in
H^{1}(G,K[V])$ contradicting the hypotheses of the proposition. Hence we
really have $m\notin
(a_{1},a_{2})_{K[V]^{G}}$. Now we show $m\Ann_{K[V]^G}(g)\subseteq
(a_{1},a_{2})_{K[V]^{G}}$. Let $a_{3}\in \Ann_{K[V]^G}(g)$. Then there is a $b_{3}\in
K[V]$ with $a_{3}g_{\sigma}=(\sigma-1)b_{3} \myforall \sigma\in G$.
Let
\begin{equation*} 
u_{ij}:=a_{i}b_{j}-a_{j}b_{i}\in K[V]^{G} \textrm{ for } 1 \le i < j \le 3.
\end{equation*}
Obviously, $m=u_{12}$ and we have
\[
u_{23}a_{1}-u_{13}a_{2}+ma_{3}=
\left|
\begin{array}{ccc}
a_{1} & a_{2} & a_{3}\\
a_{1} & a_{2} & a_{3}\\
b_{1} & b_{2} & b_{3}\\
\end{array}
\right|=0,
\]
hence $ma_{3}\in(a_{1},a_{2})_{K[V]^{G}}$. Since $a_{3}\in \Ann_{K[V]^G}(g)$ was
arbitrary, we have $m\Ann_{K[V]^G}(g)\subseteq (a_{1},a_{2})_{K[V]^{G}}$.\qed\\

Technically, the following theorem is our main result.
\begin{Satz} \label{BigMainTheorem}
Assume there is a $0 \ne g \in
H^{1}(G,K[V])$. Let  $a_{1},\ldots,a_{k}\in K[V]^{G}$ with $k\ge 2$ and
$a_{i}g=0\in H^{1}(G,K[V])$ for $i=1,\ldots,k$. 
If one of the following two conditions
\begin{enumerate}
\renewcommand{\labelenumi}{(\alph{enumi})}
\item $G$ is reductive and $a_{1},\ldots,a_{k}$ form a phsop in $K[V]$.
\item $a_{1},a_{2}$ are coprime in $K[V]$,
  and $a_{1},\ldots,a_{k}$ form a phsop in~$K[V]^{G}$.
\end{enumerate}
is true, then
\[
\cmdef K[V]^{G} \ge k-2.
\]
\end{Satz}

In the case of $(a)$ and $k=3$, this is Kemper \cite[Proposition 6]{KemperLinRed}.\\

\Bew Condition (a) implies condition (b) by  Lemma \ref{ZweiReg} and Lemma \ref{redphsop}, so let us
assume condition (b). By hypothesis,  $\height
\Ann_{K[V]^G}(g)\ge k$, and by Proposition \ref{depthKillerideal}, we have $\depth
\Ann_{K[V]^G}(g)=2$. Thus $\cmdef K[V]^{G}\ge\cmdef \Ann_{K[V]^G}(g)\ge k-2.$\qed

\section{Invariant rings with big Cohen-Macaulay defect}
In this section, given a reductive, but not linearly reductive group $G$ and a $k\in{\mathbb N}$, we will explicitly construct a
$G$-module $V$ that fulfills the hypotheses of Theorem \ref{BigMainTheorem},
hence $\cmdef K[V]^{G}\ge k-2$. The main step is the construction of a
$G$-module $U$ with a $0\ne g\in H^{1}(G,U)$. For some classical groups, this
has been done in Kohls \cite{DAPaper} by explicit calculation. With the help
of a result of Nagata, we can give a general construction.
\begin{Def} \label{Tensor}
Let $W,V$ be $G$-modules with $W\subseteq V$. Then 
\[
\Hom_{K}(V,W)_{0}:=\{ f \in \Hom_{K}(V,W): f|_{W}=0\}
\]
is a submodule of $\Hom_{K}(V,W)$, and we have
\[
\Hom_{K}(V,W)_{0} \cong W \otimes (V/W)^{*}.
\]
\end{Def}

\begin{Prop}[{Kohls \cite[Proposition 6]{DAPaper}}] \label{kompl}
Let $W$ be a submodule of a $G$-module $V$ and $\iota \in
\Hom_{K}(V,W)$ with $\iota|_{W}=\id_{W}$. Then $\sigma \mapsto
g_{\sigma}:=(\sigma-1)\iota$ is a cocycle in
$Z^{1}(G,\Hom_{K}(V,W)_{0})$, which is a coboundary if and only if
there exists a $G$-invariant complement for $W$.
\end{Prop}

Regarding this proposition, we see that in order to construct a non-trivial
cocycle, one has to find a $G$-module $V$ that contains a submodule which has
no complement. By definition, such a $G$-module $V$ exists if and only if the
group $G$ is not linearly reductive. Next, we will show how to find such a
$V$.

\begin{Def} \label{FpV}
Assume $p>0$.We call the $G$-submodule
 \[
F^{p}(V):=\left\{f\in S^{p}(V): \textrm{ there exists a
  } v\in V \textrm{ with }f=v^{p}   \right\}
\]
of $S^{p}(V)$ the \emph{$p$-th Frobenius power of $V$}.
\end{Def}

Recall that every linear algebraic group $G$ is isomorphic to a closed subgroup
of a suitable $\GL_{n}(K)$ (\cite[Theorem 2.3.6]{SpringerLin}). Then $K^{n}$
is a faithful $G$-module.

\begin{Satz}[{Nagata \cite[Proof of Theorem 1]{NagataComplete}}] \label{NagataLinRedCrit}
Let $p>0$, $G$ be a \emph{connected} linear algebraic group, and $V$ be a
\emph{faithful} $G$-module. The following are equivalent:

(a) $G$ is linearly reductive.

(b) $G$ is a torus.

(c) The submodule $F^{p}(V)$ of $S^{p}(V)$ has a complement in $S^{p}(V)$.
\end{Satz}

\begin{Korollar} \label{NagatasNoComplementCorollar}
Let $p>0$, $G$ be a linear algebraic group such that the connected
component of the unit element $G^{0}$ is not a torus, and $V$ be a \emph{faithful}
$G$-module. Then the submodule $F^{p}(V)$ of $S^{p}(V)$ has no complement in
$S^{p}(V)$. In particular, $G$ is not linearly reductive.
\end{Korollar}

Corollary \ref{NagatasNoComplementCorollar} together with Proposition
\ref{kompl} \emph{explicitly} leads to the construction of a $G$-module $U$
and a non-trivial cocycle $g\in Z^{1}(G,U)$. All we need to start is a
faithful $G$-module (which always exists). So the next step is to find annihilators of $g$. If $W$ is another
$G$-module and $w\in W^{G}$, then for a $g\in Z^{1}(G,U)$ we can define in an
obvious manner $w\otimes g\in Z^{1}(G,W\otimes U)$, and we also
get a map $w\otimes: H^{1}(G,U)\rightarrow H^{1}(G,W\otimes U)$. 

Let $g\in Z^{1}(G,V)$. Then the $K$-vector space $\tilde{V}:=V\oplus K$ can be
turned into a $G$-module with the $G$-action given by
$\sigma\cdot(v,\lambda):=(\sigma v+\lambda g_{\sigma},\lambda)$ for all
$(v,\lambda)\in \tilde{V}, \sigma\in G$. Up to $G$-module isomorphism,
$\tilde{V}$ only depends on $g+B^{1}(G,V)$. We call
$\tilde{V}$ the \emph{(corresponding) extended} $G$-module of $V$ (by $g$).

\begin{Prop}[{Kemper \cite[Proposition 2]{KemperLinRed}}] \label{anni}
Let $U$ be a $G$-module, $g \in Z^{1}(G,U)$ be a cocycle, and let
$\tilde{U}=U\oplus K$ be the extended $G$-module corresponding to $g$. Let
further be $\pi: \tilde{U}\rightarrow K, \, (u,\lambda)\mapsto \lambda$ (with
$u\in U, \lambda\in K$). Then $\pi$ is invariant, $\pi\in\tilde{U}^{*G}$, and $\pi\otimes g=0\in
H^{1}(G,\tilde{U}^{*}\otimes U)$.
\end{Prop}

\begin{Satz} \label{depthToInfty}
Let $G$ be a reductive group, $U$ be a $G$-module such that there is a $0\ne g\in H^{1}(G,U)$ and
$\tilde{U}$ be the corresponding extended $G$-module. If
$
 V:=U^{*}\oplus\bigoplus_{i=1}^{k}\tilde{U},
$
then we have
$
\cmdef K[V]^{G}\ge k-2.
$
\end{Satz}

\Bew Since $U$ is a direct summand of $K[V]=S(V^{*})$, we have (after an
embedding) $0\ne g\in H^{1}(G,K[V])$. Let  $a_{1},\ldots,a_{k}\in K[V]^{G}$ be
the $k$ copies of the element $\pi \in \tilde{U}^{*G}$ of Proposition
\ref{anni} in the $k$ summands $\tilde{U}^{*}$ of $K[V]$. Then
$a_{1},\ldots,a_{k}$ form a phsop in $K[V]$ and we have $a_{i}g=0\in
H^{1}(G,K[V])$ for $i=1,\ldots,k$. Now the result follows from Theorem \ref{BigMainTheorem}, case
(a).\qed

\begin{Satz} \label{depthInfVects}
Under the hypotheses of Theorem \ref{depthToInfty}, with $V:=U^{*}\oplus
\tilde{U}$, we have
\[
\cmdef K\left[V^{\oplus k}\right]^{G}\ge k-2 \quad \myforall k\in{\mathbb N}.
\]
\end{Satz}

If $W$ is any faithful $G$-module, then $V:=W\oplus U^{*}\oplus\tilde{U}$ is
faithful, and the theorem above remains valid with this $V$. Now let $G$ be
a reductive group that is not linearly reductive. Then by definition, there
exists a $G$-module $M$ with a submodule $N\subseteq M$ without complement. By
Proposition \ref{kompl}, the module $U:=\Hom_{K}(M,N)_{0}$ satisfies the
hypotheses of Theorem \ref{depthToInfty} or \ref{depthInfVects}, so we have proved
Theorem \ref{DasHauptErgebnis}. Together with the theorem of Hochster and
Roberts, this immediately leads to the following characterization of linearly reductive groups:

\begin{Korollar}
A reductive group $G$ is linearly reductive if and only if there is a global
  Cohen-Macaulay defect bound, i.e. a number 
  $b\in {\mathbb N}$ with $\cmdef
  K[V]^{G} \le b$ for all $G$-modules~$V$.
\end{Korollar}

Bringing all construction steps together, we
get the following explicit result. We will restrict ourselves to the case that
$G^{0}$ is not a torus. See \cite[Satz 4.2]{Doktorarbeit} for the other case.

\begin{Satz} \label{ExpliziteDarstellungen}
Let $p>0$  and $G$ be a reductive group such that $G^{0}$ is not a
  torus, and  $V$ be a faithful $G$-module. Let 
  \[U:=\Hom_{K}(S^{p}(V),F^{p}(V))_{0}\cong
  F^{p}(V)\otimes(S^{p}(V)/F^{p}(V))^{*},    \] 
(see Definition \ref{Tensor}) and $\iota \in
  \Hom_{K}(S^{p}(V),F^{p}(V))$ with $\iota|_{F^{p}(V)}=\id_{F^{p}(V)}$. 
Then with $g: G\rightarrow U, \, \sigma\mapsto (\sigma-1)\iota$ we have a
$0\ne g\in H^{1}(G,U)$. Let $\tilde{U}$ be the corresponding extended
$G$-module. Then the $G$-module 
\[
M_{k}:=F^{p}(V)^{*}\oplus(S^{p}(V)/F^{p}(V)) \oplus
  \bigoplus_{i=1}^{k}\tilde{U}
\]
is faithful, and we have
\[
\cmdef K\left[M_{k}\right]^{G} \ge k-2 \quad \myforall k\ge 0.
\]
\end{Satz}

\Bew By Proposition \ref{kompl} and Corollary
\ref{NagatasNoComplementCorollar}, we have $0\ne g\in H^{1}(G,U)$. The direct
summand $F^{p}(V)^{*}$ of $M_{k}$ is faithful since $V$ is, hence $M_{k}$ is
faithful. Since the module $U=F^{p}(V)\otimes(S^{p}(V)/F^{p}(V))^{*}$ is a direct
summand of the second symmetric power
$S^{2}(F^{p}(V)\oplus(S^{p}(V)/F^{p}(V))^{*})$, it is also a direct summand of
$K[M_{k}]=S(M_{k}^{*})$, hence after an embedding we have $0\ne g\in
H^{1}(G,K[M_{k}])$. Now the proof proceeds like the one of Theorem \ref{depthToInfty}.\qed\\

{\bf Remark.} Comparing with Theorem \ref{depthToInfty}, we see that in the definition of $M_{k}$ we replaced the summand $U^{*}=F^{p}(V)^{*}\otimes(S^{p}(V)/F^{p}(V))$ by
$F^{p}(V)^{*}\oplus(S^{p}(V)/F^{p}(V))$. This makes $M_{k}$ faithful and leads in
most cases to a lower dimension of $M_{k}$.

\section{Examples for $\SL_{2}$ and $\Ga$ invariants in positive characteristic}
The group $\SL_{2}$ acts
faithfully by left multiplication on $U:=K^{2}$. Let $\{X,Y\}$ be the standard
basis of $U$, so we have $U=\langle X,Y \rangle_{K}$. We use the notation like
$S^{2}(U)=:\langle X^{2},Y^{2},XY\rangle$ and $F^{p}(U)=:\langle
X^{p},Y^{p}\rangle$. Explicit calculations of the modules $M_{k}$ of Theorem
\ref{ExpliziteDarstellungen} lead to the following examples. These calculations
can be found in Kohls \cite[section 3]{DAPaper}, where we had to restrict
ourselves to the case $k=3$ since we did not have Theorem \ref{BigMainTheorem}.

\begin{Bsp}\label{chark2bsp}
Let $p=2$. Then
\[\cmdef K\left[ \langle X^{2},Y^{2}\rangle\oplus\bigoplus_{i=1}^{k} \langle
X^{2},Y^{2},XY\rangle  \right]^{\SL_{2}}\ge k-2\myforall k\in{\mathbb N}.\]
\end{Bsp}

\begin{Bsp}\label{chark3bsp}
Let $p=3$. Then with
\begin{eqnarray*}
&M_{k}:=&\langle X^{3},Y^{3}\rangle \oplus \langle X,Y\rangle \oplus
\bigoplus_{i=1}^{k} S^{4}(\langle X,Y\rangle)\\
\textrm{or }&M_{k}:=&\langle X^{3},Y^{3}\rangle \oplus \langle X,Y\rangle \oplus
\bigoplus_{i=1}^{k}\langle X^{2},Y^{2},XY\rangle
\end{eqnarray*}
we have $\cmdef K[M_{k}]^{\SL_{2}}\ge k-2$ for all $k\in{\mathbb N}$. In the
second case, $M_{k}$ is self-dual and completely reducible, since its summands are.
\end{Bsp}

One can use Roberts' isomorphism \cite{Roberts} to turn an example for the
group $\SL_{2}$ into an example for the additive group $\Ga=(K,+)$:  Every $\SL_{2}$-module $V$ can be regarded as a module of the
additive group $\Ga$ by the embedding $\Ga\hookrightarrow \SL_{2},
\,\,t\mapsto \left(
\begin{array}{cc}
1& t\\
0&1\\ 
\end{array}
\right)$.  Roberts' isomorphism (see \cite[Example 3.6]{Bryant}) says we then have 
\begin{equation*}
K[\langle X,Y\rangle\oplus V]^{\SL_{2}}\cong K[V]^{\Ga}.
\end{equation*}
It is probably worth remarking that in positive characteristic, it is not known if for every
$\Ga$-module $V$ the invariant ring $K[V]^{\Ga}$ is finitely generated, while
in characteristic zero this is Weitzenb{\"o}ck's Theorem \cite{Weitzenbock}.
If $V$ is a $\SL_{2}$-module and we have used Theorem~\ref{BigMainTheorem} in
case (a) to show $\cmdef K[V]^{\SL_{2}}\ge k-2$, then we also have $\cmdef
K[\langle X,Y \rangle\oplus V]^{\SL_{2}}\ge k-2$, because all the hypotheses
of the theorem made for $K[V]$ will still hold for $K[\langle X,Y
\rangle\oplus V]$. Then by Roberts' isomorphism, we also have $\cmdef
K[V]^{\Ga}\ge k-2$. In the case that $\langle X,Y\rangle$ already is a direct
summand of $V$, so $V\cong \langle X,Y
\rangle\oplus V'$ with a $\SL_{2}$-module $V'$, then Roberts' isomorphism
directly tells us $\cmdef K[V']^{\Ga}\ge k-2$. In particular, the examples \ref{chark2bsp}
and \ref{chark3bsp} for the group $\SL_{2}$ can easily be turned into examples
for the group
$\Ga$ - e.g. for $p=3$ we have $\cmdef K\left[\langle X^{3},Y^{3}\rangle\oplus
\bigoplus_{i=1}^{k}\langle X^{2},Y^{2},XY\rangle\right]^{\Ga}\ge k-2$ for all $k\ge 1$. 

If $G$ is a non-trivial, connected unipotent group, there is a surjective
algebraic homomorphism $G\rightarrow \Ga$ (see
 \cite[last paragraph before
section 3]{Bryant} for a proof of this well-known result). So if $V$ is any $\Ga$-module, it can be regarded as a
$G$-module with the same invariant ring. In particular, we have
\begin{Satz}
For every non-trivial, connected unipotent group $G$ over an algebraically
closed field $K$ of characteristic $p>0$, there exists a
$G$-module $V$ such that  $K[V^{\oplus k}]^{G}$ is finitely generated
and $\cmdef K[V^{\oplus k}]^{G}\ge k-2$ for all $k\ge 1$.
\end{Satz}

The modules $M_{k}$ of Theorem \ref{ExpliziteDarstellungen} often are not very
``nice'', in particular they have big dimensions. With some more effort, we
succeeded to construct ``nicer'' modules for
the groups $\SL_{2}$ and $\Ga$ such that the invariant rings have big
Cohen-Macaulay defect. We just state the result here, and refer to my thesis
\cite[pp. 113-126]{Doktorarbeit} for the proof.

\begin{Satz} \label{ZweitesHauptResultat}
Let $\langle X,Y \rangle$ be the natural representation of $\SL_{2}$ and $p>0$.  Let 
\[
V:=\langle X^{p},Y^{p}\rangle \oplus \bigoplus_{i=1}^{k} \langle X,Y \rangle.
\]
Then we have
\[
\cmdef K[V]^{\Ga} \ge k-2 \quad \textrm{ and } \quad \cmdef K[V]^{\SL_{2}} \ge k-3.
\]
As a direct sum of self dual $\Ga$- or $\SL_{2}$- modules, $V$ is
self-dual, too. Regarded as $\SL_{2}$-module, $V$ is completely reducible as a
direct sum of irreducible $\SL_{2}$-modules. Furthermore,
\begin{eqnarray*}
\dim K[V]^{\Ga}=2k+1 & \textrm{ and }&  \dim K[V]^{\SL_{2}}=2k-1,
\quad\textrm{ so we have}\\
\depth K[V]^{\Ga}\le k+3  & \textrm{ and } & \depth K[V]^{\SL_{2}}\le k+2.
\end{eqnarray*}
\end{Satz}

We have the conjecture that all inequations in this theorem in fact are
equations if $k\ge 2$ and $k\ge 3$ for $\Ga$ and
$\SL_{2}$ respectively. We were able to verify this conjecture with computational methods using
\Magma \cite{Magma} for
the group $\Ga$ in the 
cases  $(p,k)\in\{(2,2),(2,3),(2,4),(3,2),(3,3)\}$ - then by Roberts'
isomorphism, the conjecture is also true for $\SL_{2}$ and the corresponding
pair $(p,k+1)$. See \cite[pp. 129-140]{Doktorarbeit} for the computational
details. It is interesting to compare this theorem with the result of Hochster
that we mentioned in the text below Theorem \ref{DasHauptErgebnis}.


\begin{thebibliography}{10}

\bibitem{Magma}
Wieb Bosma, John Cannon, and Catherine Playoust.
\newblock The {M}agma algebra system. {I}. {T}he user language.
\newblock {\em J. Symbolic Comput.}, 24(3-4):235--265, 1997.
\newblock Computational algebra and number theory (London, 1993).

\bibitem{BrunsHerzog}
Winfried Bruns and J{\"u}rgen Herzog.
\newblock {\em Cohen-{M}acaulay rings}, volume~39 of {\em Cambridge Studies in
  Advanced Mathematics}.
\newblock Cambridge University Press, Cambridge, 1993.

\bibitem{Bryant}
Roger~M. Bryant and Gregor Kemper.
\newblock Global degree bounds and the transfer principle for invariants.
\newblock {\em J. Algebra}, 284(1):80--90, 2005.

\bibitem{CampEtAl}
H.~E.~A. Campbell, I.~P. Hughes, G.~Kemper, R.~J. Shank, and D.~L. Wehlau.
\newblock Depth of modular invariant rings.
\newblock {\em Transform. Groups}, 5(1):21--34, 2000.

\bibitem{Eisenbud}
David Eisenbud.
\newblock {\em Commutative algebra with a view toward algebraic geometry},
  volume 150 of {\em Graduate Texts in Mathematics}.
\newblock Springer-Verlag, New York, 1995.

\bibitem{Ellingsrud}
Geir Ellingsrud and Tor Skjelbred.
\newblock Profondeur d'anneaux d'invariants en caract\'eristique {$p$}.
\newblock {\em Compositio Math.}, 41(2):233--244, 1980.

\bibitem{FleischmannCohomology}
Peter Fleischmann, Gregor Kemper, and R.~James Shank.
\newblock On the depth of cohomology modules.
\newblock {\em Q. J. Math.}, 55(2):167--184, 2004.

\bibitem{FleischmannCohomConnectivity}
Peter Fleischmann, Gregor Kemper, and R.~James Shank.
\newblock Depth and cohomological connectivity in modular invariant theory.
\newblock {\em Trans. Amer. Math. Soc.}, 357(9):3605--3621 (electronic), 2005.

\bibitem{FleischmannDepthTrace}
Peter Fleischmann and R.~James Shank.
\newblock The relative trace ideal and the depth of modular rings of
  invariants.
\newblock {\em Arch. Math.}, 80(4):347--353, 2003.

\bibitem{Gordeev}
Nikolai Gordeev and Gregor Kemper.
\newblock On the branch locus of quotients by finite groups and the depth of
  the algebra of invariants.
\newblock {\em J. Algebra}, 268(1):22--38, 2003.

\bibitem{Hochster}
M.~Hochster.
\newblock Grassmannians and their {S}chubert subvarieties are arithmetically
  {C}ohen-{M}acaulay.
\newblock {\em J. Algebra}, 25:40--57, 1973.

\bibitem{HochsterRoberts}
Melvin Hochster and Joel~L. Roberts.
\newblock Rings of invariants of reductive groups acting on regular rings are
  {C}ohen-{M}acaulay.
\newblock {\em Advances in Math.}, 13:115--175, 1974.

\bibitem{KemperLinRed}
Gregor Kemper.
\newblock A characterization of linearly reductive groups by their invariants.
\newblock {\em Transform. Groups}, 5(1):85--92, 2000.

\bibitem{KemperDepthCohom}
Gregor Kemper.
\newblock The depth of invariant rings and cohomology.
\newblock {\em J. Algebra}, 245(2):463--531, 2001.
\newblock With an appendix by Kay Magaard.

\bibitem{DAPaper}
M.~Kohls.
\newblock Non {C}ohen-{M}acaulay invariant rings of infinite groups.
\newblock {\em J. Algebra}, 306:591--609, 2006.

\bibitem{Doktorarbeit}
M.~Kohls.
\newblock {\"U}ber die {T}iefe von {I}nvariantenringen unendlicher {G}ruppen.
\newblock {\em Doktorarbeit, Technische Universit\"at M\"unchen. Available
  from\\
  {http://nbn-resolving.de/urn/resolver.pl?urn:nbn:de:bvb:91-diss-20071116-619%
331-1-3} or {arXiv:0711.4730v1}}, pages 1--148, 2007.

\bibitem{NagataComplete}
Masayoshi Nagata.
\newblock Complete reducibility of rational representations of a matric group.
\newblock {\em J. Math. Kyoto Univ.}, 1:87--99, 1961/1962.

\bibitem{NagataAffine}
Masayoshi Nagata.
\newblock Invariants of a group in an affine ring.
\newblock {\em J. Math. Kyoto Univ.}, 3:369--377, 1963/1964.

\bibitem{NeuselDepth}
M.~D. Neusel.
\newblock Comparing the depths of rings of invariants.
\newblock In {\em Invariant theory in all characteristics}, volume~35 of {\em
  CRM Proc. Lecture Notes}, pages 189--192. Amer. Math. Soc., Providence, RI,
  2004.

\bibitem{Roberts}
Michael Roberts.
\newblock On the {C}ovariants of a {B}inary {Q}uantic of the $n^{th}$ {D}egree.
\newblock {\em The Quarterly Journal of Pure and Applied Mathematics},
  4:168--178, 1861.

\bibitem{ShankWehlau}
R.~James Shank and David~L. Wehlau.
\newblock On the depth of the invariants of the symmetric power representations
  of {${\rm SL}\sb 2({\bf F}\sb p)$}.
\newblock {\em J. Algebra}, 218(2):642--653, 1999.

\bibitem{SpringerLin}
T.~A. Springer.
\newblock {\em Linear algebraic groups}, volume~9 of {\em Progress in
  Mathematics}.
\newblock Birkh\"auser Boston, Mass., 1981.

\bibitem{Weitzenbock}
R.~Weitzenb\"ock.
\newblock {\"U}ber die {I}nvarianten von linearen {G}ruppen.
\newblock {\em Acta. Math.}, 58:231--293, 1932.

\end{thebibliography}
\end{document}